\documentclass{article}
\usepackage{amsfonts}
\usepackage{amsmath,amssymb,pxfonts, color}

\newtheorem{lem}{Lemma}[section]
\newtheorem{proposition}{Proposition}[section]
\newtheorem{thm}{Theorem}[section]
\newtheorem{Rem}{Remark}[section]
\begin{document}

\title{The asymptotic expansion of the regular discretization error of
 It\^o integrals}

\author{Elisa Al\`{o}s\thanks{%
Supported by the Spanish grant MTM2013-40782-P.} \\
Dpt. d'Economia i Empresa\\
Universitat Pompeu Fabra \\
and Barcelona GSE\\
c/Ramon Trias Fargas, 25-27\\
08005 Barcelona, Spain \\
elisa.alos@upf.edu
\and
Masaaki Fukasawa \\
Graduate School of Engineering Science\\
Osaka University\\
1-3 Machikaneyama, Toyonaka\\
560-8531 Osaka, Japan\\
fukasawa@sigmath.es.osaka-u.ac.jp}

\maketitle

\begin{abstract}
We study a Edgeworth-type refinement of the central limit theorem for the discretizacion error of It\^{o} integrals. Towards this end, we introduce a new approach,  based on the anticipating It\^{o} formula. This alternative technique allows us to compute explicitly the terms of the corresponding expansion formula.

\vspace{2mm}


\vspace{2mm}


\end{abstract}
\section{Introduction}
The It\^o integral for semi-martingales is defined as a limit of a random sequence of
Riemann-Stieltjes type sum. First it is defined as an $L^2$ limit and
then, the limit is characterized also in terms of convergence in
probability in a functional sense.
A natural question, from both theoretical and practical viewpoints,
is how close the limit and the Riemann-Stieltjes type
sum are. The difference between them, which we call the discretization error,
should converge to 0 as the time partition for the Riemann-Stieltjes sum
becomes finer and finer. The question is  how fast  it is  and how a
renormalized error behaves.

The first answer to this question was given by Rootz\'en~\cite{Rootzen},
where a central limit theorem for the regular discretization error, that
is, the discretization error when the time partition is regular; $0,
1/n, 2/n, \dots $ was proven. The convergence rate is $\sqrt{n}$ and the
limit distribution is mixed normal.
The convergence holds stably and the theory of stable convergences for
discretized processes to conditionally Gaussian processes was further
developed by Jacod and Shiryaev~\cite{JS}.
The theory has remarkable applications to high-frequency data analysis;
see e.g. A\"it-Sahalia and Jacod~\cite{AJ} and Jacod and
Protter~\cite{JP}.
The general theory also enables us to treat irregular and stochastic
time partitions for the Riemann-Stieltjes sum; see
Fukasawa~\cite{F2011D}.
In the context of mathematical finance, the discretization error
corresponds to the discrete hedging error, which is inevitably
associated with the discretization of hedging strategies.
The limit distribution tells us how to quantify the risk of hedging by
finite transactions that is usually negleted in continuous-time
financial modeling.

There are two directions for a refinement of the central limit theorem;
the large deviation theory and the theory of the Edgeworth expansion.
The former focuses the tail probability and gives a precise asymptotic
formula in terms of the so-called rate function.
The latter focuses the behavior around the mean and gives an expansion
formula in terms of Hermite polynomials with coefficients determined by
moments.
The former tends to give a more precise approximation, while it requires
to compute the rate function that is impossible in many problems.
The latter is often less precise; however it only requires to compute
moments that is usually feasible. Further, in the context of statistics,
the Edgeworth expansion provides a theoretical justification of some
popular computational techniques like the bootstrap; see
Hall~\cite{Hall}.
Unfortunaly, 
both of the refinements are not directly applicable to the discrezation
error of It\^o integral  because the limit distribution is in general not normal but only
mixed normal.

The theory of the Edgeworth expansion in mixed normal limit has
recently developed by Yoshida~\cite{YoshidaMixed} by extending his
martingale expansion approach~\cite{YoshidaMart1,YoshidaMart2} in normal limit.
It turns out that the expansion formula involves entangled anticipating
effects described in terms of the Malliavin calculus.
The Yoshida theory provides a general framework and for each of concrete applications, a
non-negligble effort is still required to obtain an explicit expression
of the expansion terms.
An application to the power variation of diffusion processes is given by
Podolskij and Yoshida~\cite{PY}.

The aim of this study is to give an explicit expansion formula for the
regular discretization error of It\^o integrals.
We treat a  far more restricted problem than in
Yoshida~\cite{YoshidaMixed} and try to get a more explicit formula.
For this purpose, instead of just analyzing abstract terms given in
Yoshida~\cite{YoshidaMixed}, we introduce a more systematic approach
based on the anticipating It\^o formula, developed by
Nualart and Pardoux in~\cite{NP}. This alternative approach is in a sense more
elementary and the appearance of Hermite polynomials is more natural.

\section{Anticipating stochastic calculus}

In this section we recall the basic results on the Malliavin calculus and
the anticipating stochastic calculus
 we use
through the paper. For a more detailed introduction to this subject we
refer to \cite{Nualart}.

 Let $T>0$ and $(\Omega, \mathcal{F}, P, \{\mathcal{F}_t\}_{0
 \leq t \leq T})$ be an filtered probability
space satisfying the usual conditions which supports a standard Brownian
motion $W$ on $[0,T]$. We denote by $E$ the expectation operator with respect to $P$,
and denote by
$D$  the Malliavin derivative operator with respect to $W$.
More precisely, we will assume a partial Malliavin structure
$(\Omega,\mathcal{F},P) = (\Omega_1\times \Omega_2, \mathcal{F}_1
 \otimes \mathcal{F}_2, P_1 \otimes P_2)$
 with $(\Omega_1,\mathcal{F}_1,P_1)$ being the Wiener space associated with $W$.
It is
 well-known that $D$ is a  closable operator from $L^p(\Omega)$ to
 $L^p([0,T]\times\Omega)$, for any $p\geq 1$. We will denote by
 $\mathbb{D}^{1,p}$ the domain of $D$ in $L^p(\Omega)$  . We also consider
 the iterated derivatives $D^{n}$, for $n\ge 1$, whose domains will be
 denoted by $\mathbb{D}^{n,p}$.
The Sobolev norm of $\mathbb{D}^{n,p}$ will be denoted by $\|\cdot\|_{n,p}$.
 We will use the notation $\mathbb{L}^{n,p}= L^2([0,T];\mathbb{D}^{n,p})$.

 Given a process $X\in \mathbb{L}^{1,p}$, $D^+ X$ and $D^-X$ will be the
 element of $L^p([0,T]\times\Omega)$ satisfying 
 \begin{equation}
  \label{dertt1}
\lim_{n \to \infty} \int_0^T \sup_{s < t \leq (s + \frac{1}{n})\wedge T}
E[|D_sX_t - D^+  X_t|^p] \mathrm{d}t= 0
 \end{equation}
 and
 \begin{equation}
    \label{dertt2}
\lim_{n \to \infty} \int_0^T \sup_{(s-\frac{1}{n})\vee 0 \leq  t < s}
E[|D_sX_t - D^-  X_t|^p] \mathrm{d}t= 0
 \end{equation}
respectively.
Moreover, we will denote by $\mathbb{L}_{+}^{1,p}$ and
$\mathbb{L}_{-}^{1,p}$ the sets of such processes $X$ in $\mathbb{L}^{1,p}$
that admit $D^+X$ and $D^-X$ with
(\ref{dertt1}) and (\ref{dertt2}) respectively.

The adjoint of the derivative operator $D$, denoted by
$\delta$, is an extension of the It\^{o} integral in the sense that
the set $L_{a}^{2}([0,T]\times \Omega )$ of square integrable and
$\{\mathcal{F}_t\}$-adapted processes is included in the domain of $\delta$ and the operator $\delta$
restricted to $L_{a}^{2}([0,T]\times \Omega )$ coincides
with the It\^{o} integral. We will make use of the notation
\begin{equation*}
\delta(u)=\int_{0}^{T}u_{t}\mathrm{d}W_{t}.
\end{equation*}
We recall that
$\mathbb{L}^{n,2}$ is included
in the domain of $\delta$ for all $n\geq 1$.

\begin{proposition} (The Malliavin derivative of an It\^{o} process)
\label{MalliavinderItoprocess}
Consider an It\^{o} process of the form
$$
S_{t}=S_{0}+\int_{0}^{t}S_{u}^{\prime }\mathrm{d}u+\int_{0}^{t}S_{u}^{%
\prime \prime }\mathrm{d}W_{u}, 
$$
where $S_{0}$ is a positive constant and $S^{\prime }, S^{\prime \prime
 }$ are adapted processes in $\mathbb{L}^{1,2}$. Then
 $S\in\mathbb{L}^{1,2}$ and for all $0< r<t<T$, 
$$
D_rS_{t}=\int_{r}^{t}D_rS_{u}^{\prime }\mathrm{d}u+S_{r}^{\prime \prime }+\int_{r}^{t}D_rS_{u}^{\prime \prime }\mathrm{d}W_{u}.
$$
\end{proposition}

The key tool for this work is the following anticipating It\^{o}
formula;
see Theorem 3.2.4 of \cite{Nualart}.   

\begin{thm}
 \label{ito}
 Let $S$ and $A$ be processes of the form 
\begin{equation*}
\begin{split}
& S_{t}=S_{0}+\int_{0}^{t}S_{u}^{\prime }\mathrm{d}u+\int_{0}^{t}S_{u}^{%
\prime \prime }\mathrm{d}W_{u}, \\
& A_{t}=\int_{t}^{T}A_{s}^{\prime }\mathrm{d}s,
\end{split}%
\end{equation*}%
where $S^{\prime },S^{\prime \prime }\in L_{a}^{2}([0,T]\times \Omega )$ and 
$A^{\prime }\in \mathbb{L}^{1,2}$. Then, $A\in \mathbb{L}_{-}^{1,2}$,
 $$D^{-}A_{u}=\int_{u}^{T}D_{u}A_{s}^{\prime }\mathrm{d}s$$ and
 for any $f\in C_{b}^{2}(\mathbb{R}^{2})$ with
\begin{equation*}
E\left[\left( \int_{0}^{T}|f_{1,1}(A_{u},S_{u})D^{-}A_{u}|\mathrm{d}\langle
S,W\rangle _{u}\right)^{2} \right] <\infty,
\end{equation*}
 it holds
\begin{equation*}
\begin{split}
f(A_{t},S_{t})=& f(A_{0},S_{0})+\int_{0}^{t}f_{1,0}(A_{u},S_{u})\mathrm{d}%
A_{u}+\int_{0}^{t}f_{0,1}(A_{u},S_{u})\mathrm{d}S_{u} \\
& +\int_{0}^{t}f_{1,1}(A_{u},S_{u})D^{-}A_{u}\mathrm{d}\langle S,W\rangle
_{u}+\frac{1}{2}\int_{0}^{t}f_{0,2}(A_{u},S_{u})\mathrm{d}\langle S\rangle
_{u}
\end{split}%
\end{equation*}%
 for all $t\in \lbrack 0,T]$, where 
\begin{equation*}
\int_{0}^{t}f_{1,0}(A_{u},S_{u})\mathrm{d}S_{u}=%
\int_{0}^{t}f_{1,0}(A_{u},S_{u})S_{u}^{\prime }\mathrm{d}u+\delta
(1_{[0,t]}f_{1,0}(A,S)S^{\prime \prime }).
\end{equation*}%

\end{thm}

\section{The statement of the main result}
We consider the limit distribution of a normalized regular
discretization error of a stochastic integral 
\begin{equation*}
Z^n_t = \sqrt{n} \left\{ \int_0^t X_s \mathrm{d}Y_s - \sum_{j=0}^{\infty}
X_{t^n_j}(Y_{t^n_{j+1}\wedge t} - Y_{t^n_j\wedge t}) \right\}
\end{equation*}
as $n\to \infty$, where $t^n_j = j/n$ and $X$ and $Y$ are continuous It\^o
processes of the form 
\begin{equation}
\label{sde1}
\mathrm{d}X_t = \Xi_t \mathrm{d}t + \Gamma_t\mathrm{d}W_t, \ \ \mathrm{d}Y_t
= \Theta_t \mathrm{d}t + \Sigma_t\mathrm{d}W_t, \ \ 
\end{equation}
where $\Xi, \Theta, \Gamma, \Sigma$ are continuous processes adapted to the
filtration $\{\mathcal{F}_t\}$. We further assume that the product $\Gamma
\Sigma$ is not identically zero and that $\Gamma$ and $\Sigma^2$ (the square
of $\Sigma$) are also It\^o processes of the form 
\begin{equation}
\label{eq}
\mathrm{d}\Gamma_t = \check{\Gamma}_t\mathrm{d}t + \hat{\Gamma}_t\mathrm{d}%
W_t, \ \ \mathrm{d}\Sigma^2_t = \check{\Sigma}_t\mathrm{d}t + \hat{\Sigma}_t%
\mathrm{d}W_t,
\end{equation}
where $\check{\Gamma}, \check{\Sigma}, \hat{\Gamma}$ and $\hat{\Sigma}$ are
adapted processes.
\vspace{0.1cm}

We will assume the following hypotheses:
\begin{description}
\item[(H1)] For all $p>0,$ the adapted processes $\Xi ,\Theta ,\check{\Gamma}%
	   ,\hat{\Gamma},\check{\Sigma},\hat{\Sigma}\ $are in $\mathbb{L}^{4,p}.$
\item[(H2)]
	   $\Theta$ and $\Sigma$ belong to $\mathbb{L}^{1,p}_{+}$ for
	   all $p>2$
	   and $D^+\Theta$ and $D^+\Sigma$ are continuous.

	   \item[(H3)] For any $U \in \{\Xi, \Theta, \hat{\Sigma}, \hat{\Gamma} \}$, 
\begin{equation*}
\sup_{t \in [0,T]} E[|U_t|^{4}] < \infty, \ \ \lim_{n \to
\infty}\sup_{|s-t|\leq \frac{1}{n}} E[|U_t - U_s|^2] = 0.
\end{equation*}
\end{description}

\begin{Rem}
All of the above assumptions are satisfied, for example, if $X_t =
g_X(t,W_t) $ and $Y_t = g_Y(t,W_t)$ for some $C^\infty$ functions $g_X,
 g_Y$ with all the derivatives being of at most exponential growth. This is the case of discrete hedging
for a European $C^\infty$ payoff of at most polynomial growth under the Black-Scholes model.
\end{Rem}
\begin{Rem}
\label{lem4p} Under (H1), $\Gamma $, $\Sigma ^{2}$,  $X$ and $Y$ are adapted processes
in $\mathbb{L}^{4,p}$ for all $p>0$. Moreover, for all $p>0$,
$$\sup_{t \in
[0,T]}E[|\Gamma_t|^p] + \sup_{t \in [0,T]}E[|\Sigma_t|^p] <\infty.$$
\end{Rem}

We will also make use of the following notation
$$
V_t = \frac{1}{2}\int_t^T \Gamma_s^2\Sigma_s^2\mathrm{d}s.
$$
Using Proposition \ref{MalliavinderItoprocess}, it is easy to see that
$V\in\mathbb{L}_{-}^{1,p}$ for any $p>1$ and 
$$
D^-V_t = \frac{1}{2}\int_t^T D_t\left[\Gamma_s^2\Sigma_s^2\right]\mathrm{d}s.
$$
Moreover,  $D^-V$ and  $D\left[\Gamma_s^2\Sigma_s^2\right]$ $\in \mathbb{L}_-^{1,p}$ for any $s$.

\vspace{2mm}

 Now we can state the main result of this paper.
\begin{thm}
\label{main} Under (H1), (H2) and (H3), for any $f \in C^\infty_b(\mathbb{R}%
) $, we have 
\begin{equation*}
E[f(Z^n_T)] = \int_\mathbb{R}f(z)E[Q_n(z)]\mathrm{d}z + o(n^{-1/2}),
\end{equation*}
where 
\begin{equation*}
Q_n(z) = \left\{1 + \frac{1}{\sqrt{n}}\left\{ A_1 H_1(z,V_0) + A_3H_3(z,V_0)
+ A_5H_5(z,V_0) \right\}\right\}\phi(z,V_0),
\end{equation*}
\begin{equation*}
\begin{split}
& H_k(z, t) = \frac{(-1)^k}{\phi(z, t)} \frac{\partial^k}{\partial z^k}
\phi(z, t), \ \ k=1,3,5, \\
&\phi(z, t) = \frac{1}{\sqrt{2\pi t}} \exp\left\{ -\frac{z^2}{2 t}
\right\}
\end{split}%
\end{equation*}
and 
\begin{equation*}
\begin{split}
&A_1 = \frac{1}{2}\int_0^T (\Xi_t\Theta_t + D^+\Theta_t\Gamma_t)\mathrm{d}t,
\\
& A_3 = \frac{1}{4} \int_0^T \left( (\Xi_t
 \Sigma_t + \Theta_t \Gamma_t + 
D^+\Sigma_t\Gamma_t)D^-V_t + 
\Gamma_t \Sigma_t  (D^-)^2V_t \right) \mathrm{d%
}t +\frac{1}{6}\int_0^T \Sigma_t^3\Gamma_t^3 \mathrm{d}t,\\
&A_5 = \frac{1}{8}\int_0^T  \Gamma_t\Sigma_t |D^-V_t|^2 \mathrm{d}t.
\end{split}%
\end{equation*}
\end{thm}
This is an Edgeworth type expansion in that the expansion coefficients are
written in terms of the normal density and Hermite polynomials. An
important difference from the classical Edgeworth expansion formula  is that the limit variance $V_0$ is random.
Due to this anticipating feature, the 5th order Hermite polynomial
appears in this first order expansion, while
Hermite polynomials of only up to 3rd order appear  
 in the first order Edgeworth expansion in the classical situation, that is,
 where $V_0$ is deterministic.
The 1st and 3rd order Hermite polynomial terms represent, respectively,
 the bias and the skewness of the limit distribution of $O(n^{-1/2})$.

\section{The outline of the proof}

\label{sec:main}

Let $X^n_t = X_{[nt]/n}$ and 
\begin{equation*}
V^n_t = n \int_t^T(X_s-X^n_s)^2\Sigma_s^2\mathrm{d}s.
\end{equation*}
Then we have 
\begin{equation}  \label{zeq}
Z^n_t = \sqrt{n}\int_0^t (X_s-X^n_s)\mathrm{d}Y_s, \ \ \mathrm{d}\langle Z^n
\rangle_t = - \mathrm{d}V^n_t.
\end{equation}
For $f \in C^\infty_b(\mathbb{R})$, define 
\begin{equation*}
q(t,x) = \int f(z)\phi(x-z,t)\mathrm{d}z.
\end{equation*}
Note that $q$ is the solution of the heat equation 
\begin{equation}  \label{heq}
q_{1,0} = \frac{1}{2}q_{0,2}, \ \ q(0,x) = f(x),
\end{equation}
which in particular implies that $q \in C^\infty_b([0,\infty)\times \mathbb{R%
})$.

\vspace{0.3cm}

By (\ref{zeq}), (\ref{heq}),(H1), (H2), (H3) and Theorem~\ref{ito},%
\begin{equation*}
\begin{split}
E[f(Z_{T}^{n})]& =E[q(V_{T}^n,Z_{T}^{n})] \\
& =E[q(V_{0}^{n},0)]+\sqrt{n}%
\int_{0}^{T}E[q_{0,1}(V_{t}^{n},Z_{t}^{n})(X_{t}-X_{t}^{n})\Theta _{t}%
\mathrm{d}t] \\
& \hspace*{2cm}+\sqrt{n}\int_{0}^{T}E[q_{1,1}(V_{t}^{n},Z_{t}^{n})D^-V_{t}^{n}
(X_{t}-X_{t}^{n})\Sigma _{t}]\mathrm{d}t,
\end{split}%
\end{equation*}%
where 
\begin{equation*}
D^- V_{t}^{n}=n\int_{t}^{T}D_{t}\left[ (X_{s}-X_{s}^{n})^{2}\Sigma _{s}^{2}%
\right] \mathrm{d}s.
\end{equation*}%
We will show in Lemma ~\ref{lemT1} that 
\begin{equation*}
E[q(V_{0}^{n},0)]=E[q(V_{0},0)]+o(n^{-1/2}),
\end{equation*}%
in Lemma~\ref{lem2} that 
\begin{equation*}
\begin{split}
& n\int_{0}^{T}E[q_{0,1}(V_{t}^{n},Z_{t}^{n})(X_{t}-X_{t}^{n})\Theta _{t}]%
\mathrm{d}t \\
& \rightarrow \frac{1}{2}\int_{0}^{T}E[q_{0,1}(V_{0},0)(\Xi _{t}\Theta
_{t}+D^{+}\Theta _{t}\Gamma _{t})]\mathrm{d}t+\frac{1}{2}%
\int_{0}^{T}E[q_{1,1}(V_{0},0)D^-V_{t}\Theta _{t}\Gamma _{t}]\mathrm{d}t \\
& =E[q_{0,1}(V_{0},0)A_{1}]+\frac{1}{4}%
E\left[q_{0,3}(V_{0},0)\int_{0}^{T}D^-V_{t}\Theta _{t}\Gamma _{t}\right]\mathrm{d}t
\end{split}%
\end{equation*}%
and in Lemma \ref{lem3} that 
\begin{equation*}
\begin{split}
&
n\int_{0}^{T}E[q_{1,1}(V_{t}^{n},Z_{t}^{n})D^{-}V_{t}^{n}(X_{t}-X_{t}^{n})%
\Sigma _{t}]\mathrm{d}t \\
& \rightarrow \frac{1}{2}E\left[q_{1,1}(V_{0},0)\int_{0}^{T}(D^{-}V_{t}\Xi
_{t}\Sigma _{t}+(D^{-})^2V_{t}\Sigma _{t}\Gamma _{t}+ D^{-}V_{t}D^{+}\Sigma
_{t}\Gamma _{t})\mathrm{d}t\right] \\
& \hspace*{2cm}+\frac{1}{2}E\left[q_{2,1}(V_{0},0)\int_{0}^{T}|D^{-}V_{t}|^{2}%
\Gamma _{t}\Sigma _{t}]\mathrm{d}t \right]+  \frac{1}{3} E\left[q_{1,1} (V_0,0)\int_0^T  \Sigma_t^3\Gamma_t^3 dt\right]               \\
& =\frac{1}{4}E\left[q_{0,3}(V_{0},0)\int_{0}^{T}(D^{-}V_{t}\Xi _{t}\Sigma
_{t}+(D^{-})^2V_{t}\Sigma _{t}\Gamma _{t}+D_{t}^{-}V_{t}D^{+}\Sigma _{t}\Gamma
_{t}+\frac{2}{3} \Sigma_t^3\Gamma_t^3)\mathrm{d}t\right]\\
& \hspace*{2cm}+E[q_{0,5}(V_{0},0)A_{5}].
\end{split}%
\end{equation*}%
Since 
\begin{equation*}
q_{0,k}(t,x)=\int_{\mathbb{R}}f(z)\frac{\partial ^{k}}{\partial x^{k}}\phi
(x-z,t)\mathrm{d}z=\int_{\mathbb{R}}f(z)H_{k}(z-x,t)\phi (z-x,t)\mathrm{d}z
\end{equation*}%
for $k=0,1,2,\dots $, the expansion claimed in Theorem \ref{main} follows.%
\newline

\section{The approximating processes}
This section is devoted to prove some results of the approximating processes  $V^n$ and $Z^n$ that we will use through the paper.
We recall the following lemma, which will be repeatedly used in the sequel:
\begin{lem}
\label{Hol}  Let $p\geq 1$, $p^\prime = p/(p-1)$, $F \in
L^p([0,T]\times\Omega)$ and $G \in L^{p^\prime}([0,T]\times \Omega)$. Denote
their norms as 
\begin{equation*}
\|F\|_p = \left\{ \int_0^T E[|F_t|^p]\mathrm{d}t \right\}^{1/p}, \ \
\|G\|_{p^\prime} = \left\{ \int_0^T E[|G_t|^{p\prime}]\mathrm{d}t
\right\}^{1/p^\prime}.
\end{equation*}
Then, for all $n \in \mathbb{N}$, 
\begin{equation*}
E\left[ n \int_0^T \int_{[nt]/n}^{t}|F_sG_t| \mathrm{d}s\mathrm{d}t \right]
\leq \|F\|_p \|G\|_{p^\prime}.
\end{equation*}
\end{lem}
\textit{Proof: } By H\"older's inequality, 
\begin{equation*}
E\left[ n \int_0^T \int_{[nt]/n}^{t}|F_sG_t| \mathrm{d}s\mathrm{d}t \right]
\leq \|G\|_{p^\prime} E\left[\int_0^T \left|n\int_{[nt]/n}^t|F_s| \mathrm{d}%
s\right|^p\mathrm{d}t\right]^{1/p}.
\end{equation*}
Since $n(t-[nt]/n)\leq 1$, by Jensen's inequality, 
\begin{equation*}
E\left[\int_0^T \left|n\int_{[nt]/n}^t|F_s| \mathrm{d}s\right|^p\mathrm{d}t%
\right] \leq E\left[\int_0^T n \int_{[nt]/n}^t|F_s|^p \mathrm{d}s\mathrm{d}t%
\right].
\end{equation*}
Since $[a]\leq b$ if and only if $a < [b]+1$ in general for $a,b \geq 0$, $%
[nt]/n \leq s \leq t$ is equivalent to $s \leq t < ([ns]+1)/n$. Therefore by
the Fubini theorem, 
\begin{equation*}
E\left[\int_0^T n\int_{[nt]/n}^t|F_s|^p \mathrm{d}s\mathrm{d}t\right] = E%
\left[\int_0^T (1+[ns]-ns) |F_s|^p \mathrm{d}s\right] \leq \|F\|_p^p.
\end{equation*}
\hfill////

\begin{lem}
\label{conti} Under (H1) and (H3)
 \begin{equation*}
  \sup_{ | s-s^\prime|< \delta }\|\Sigma^2_s
-  \Sigma_{s^\prime}^2\|_{1,4} \to 0
 \end{equation*}
and
\begin{equation*}
 \sup_{ | s-s^\prime|< \delta }\|\Sigma^2_s\Gamma^2_s
-  \Sigma_{s^\prime}^2\Gamma_{s^\prime}^2\|_{1,4} \to 0
\end{equation*}
as $\delta \to 0$
\end{lem}
\textit{Proof:}
For the sake of simplicity we can assume that $s>s^\prime$. Using (H1) we can see that
\begin{equation*}
D_t\Sigma^2_s- D_t \Sigma_{s^\prime}^2=\int_{s'}^s D_t \check\Sigma_r
+\hat\Sigma_t\textbf{1}_{[s^\prime , s]}(t)
+\int_{s'}^s D_t \hat\Sigma_r dW_r.
\end{equation*}%
Then, (H3) allows us to show the first convergence. The second
convergence is treated similarly.  \hfill////

\subsection{The stable convergence of $V^n$}

The goal of this subsection is to prove  Proposition~\ref{stableconv} below, which is
necessary to prove Lemmas ~\ref{lemT1} and  \ref{lem2} in Section~\ref{leading}.

\begin{lem}
\label{lem0} Let us consider two processes $\Sigma$ and $\Gamma$ defined as in (\ref{eq}). Then, under (H1),
\begin{equation*}
\sqrt{n}\left( \int_0^T(\Gamma _{t}^{2}\Sigma _{t}^{2}-\Gamma
_{[nt]/n}^{2}\Sigma _{[nt]/n}^{2})\mathrm{d}t\right)
\end{equation*}%
converges to $0$ in $L^1$ as $n\rightarrow \infty $.
\end{lem}
\textit{Proof : } By the assumption on $\Gamma $ and $\Sigma $, $%
\Gamma^{2}\Sigma ^{2}$ is an It\^{o} semimartingale and so, of the form 
\begin{equation*}
\mathrm{d}(\Gamma ^{2}\Sigma ^{2})_{t}=\alpha _{t}\mathrm{d}t+\beta _{t}%
\mathrm{d}W_{t}
\end{equation*}%
with $\alpha, \beta \in L^2_a([0,T]\times \Omega)$. By Lemma~\ref{Hol}, 
\begin{equation*}
\sqrt{n} \int_0^T\int_{[nt]/n}^{t}\alpha _{s}\mathrm{d}s\mathrm{d}t \to 0
\end{equation*}
in $L^1$. Further, denoting $t_j = t^n_j \wedge T$, 
\begin{equation*}
\begin{split}
nE\left[ \left\vert \sum_j\int_{t_{j}}^{t_{j+1}}\int_{t_{j}}^{t}\beta _{s}%
\mathrm{d}W_{s}\mathrm{d}t\right\vert ^{2}\right] & =n\sum_jE\left[
\left\vert \int_{t_{j}}^{t_{j+1}}\int_{t_{j}}^{t}\beta _{s}\mathrm{d}W_{s}%
\mathrm{d}t\right\vert ^{2}\right] \\
& =n\sum_jE\left[ \int_{t_{j}}^{t_{j+1}}\int_{t_{j}}^{t_{j+1}}%
\int_{t_{j}}^{t}\beta _{s}\mathrm{d}W_{s}\int_{t_{j}}^{u}\beta _{v}\mathrm{d}%
W_{v}\mathrm{d}u\mathrm{d}t\right] \\
& =2n\sum_j\int_{t_{j}}^{t_{j+1}}\int_{t_{j}}^{t}\int_{t_{j}}^{u}E[ \beta
_{v}^{2}]\mathrm{d}v\mathrm{d}u\mathrm{d}t \\
&= 2n\int_0^T\int_{[nt]/n}^{t}\int_{[nt]/n}^uE[\beta_v^2]\mathrm{d}v\mathrm{d%
}u\mathrm{d}t \\
&=2n \int_0^T\int_{[nt]/n}^{t}(t-v)E[\beta_v^2]\mathrm{d}v\mathrm{d}t \\
&\leq  2 \int_0^T\int_{[nt]/n}^{t}E[\beta_v^2]\mathrm{d}v\mathrm{d}t \to 0
\end{split}%
\end{equation*}%
by  Lemma~\ref{Hol}, which completes the proof. \hfill////

\begin{proposition}
\label{stableconv}  Assume that hypotheses (H1) and (H3) hold. Then $\ \sqrt{n}(V_{0}^{n}-V_{0})$ is uniformly integrable and 
\begin{equation*}
\sqrt{n}(V_{0}^{n}-V_{0})\rightarrow \mathcal{MN}\left( 0,\frac{1}{3}%
\int_{0}^{T}\Gamma _{t}^{4}\Sigma _{t}^{4}\mathrm{d}t\right)
\end{equation*}
stably as $n\rightarrow \infty $. In particular, for any uniformly bounded
random variables $U_n$, $U$ such that $U_n \to U$ in probability, 
\begin{equation*}
E[\sqrt{n}(V_{0}^{n}-V_{0})U_n] \to 0.
\end{equation*}
\end{proposition}
\textit{Proof : } We will denote $t_j = t^n_j \wedge T$ for brevity. Denote
by $E_{t_j}$ the conditional expectation given $\mathcal{F}_{t_j}$. We will
use $C$ as a generic constant. 
\begin{equation*}
\begin{split}
V_{0}^{n}-V_{0}=& n\int_{0}^{T}(X_{t}-X_{t}^{n})^{2}\Sigma_{t}^{2}\mathrm{d}%
t-\frac{1}{2}\int_{0}^{T}\Gamma_{t}^{2}\Sigma _{t}^{2}\mathrm{d}t \\
=& n\sum_j\int_{t_{j}}^{t_{j+1}}\left( \int_{t_j}^t\Xi_u\mathrm{d}u\right)
^{2}\Sigma _{t}^{2}\mathrm{d}t + 2n\sum_j\int_{t_{j}}^{t_{j+1}}
\int_{t_j}^t\Xi_u\mathrm{d}u \int_{t_j}^t\Gamma_u\mathrm{d}W_u\Sigma _{t}^{2}%
\mathrm{d}t \\
&+ n\sum_j\int_{t_{j}}^{t_{j+1}}\left( \Gamma_{t_{j}}\int_{t_{j}}^{t}\mathrm{%
d}W_{u}+\int_{t_{j}}^{t}(\Gamma_{u}-\Gamma _{t_{j}})\mathrm{d}W_{u}\right)
^{2}\Sigma _{t}^{2}\mathrm{d}t -\frac{1}{2}\int_{0}^{1}\Gamma_{t}^{2}\Sigma
_{t}^{2}\mathrm{d}t
\end{split}%
\end{equation*}
The sum of the last two terms can be written as 
\begin{equation*}
\begin{split}
&n\sum_j\Gamma _{t_{j}}^{2}\Sigma
_{t_{j}}^{2}\int_{t_{j}}^{t_{j+1}}(W_{t}-W_{t_{j}})^{2}\mathrm{d}t-\frac{1}{2%
}\int_{0}^{1}\Gamma_{t}^{2}\Sigma_{t}^{2}\mathrm{d}t \\
&
+n\sum_j\Gamma_{t_{j}}^{2}\int_{t_{j}}^{t_{j+1}}(W_{t}-W_{t_{j}})^{2}(%
\Sigma_{t}^{2}-\Sigma _{t_{j}}^{2})\mathrm{d}t \\
& +2n\sum_j\Gamma
_{t_{j}}\int_{t_{j}}^{t_{j+1}}(W_{t}-W_{t_{j}})\int_{t_{j}}^{t}(\Gamma
_{u}-\Gamma _{t_{j}})\mathrm{d}W_{u}\Sigma _{t}^{2}\mathrm{d}t \\
& +n\sum_j\int_{t_{j}}^{t_{j+1}}\left( \int_{t_{j}}^{t}(\Gamma _{u}-\Gamma
_{t_{j}})\mathrm{d}W_{u}\right) ^{2}\Sigma _{t}^{2}\mathrm{d}t.
\end{split}%
\end{equation*}

\noindent Step 1). First we show 
\begin{equation*}
\begin{split}
\sqrt{n}& \left( n\sum_j\Gamma _{t_{j}}^{2}\Sigma
_{t_{j}}^{2}\int_{t_{j}}^{t_{j+1}}(W_{t}-W_{t_{j}})^{2}\mathrm{d}t-\frac{1}{2%
}\int_{0}^{1}\Gamma _{t}^{2}\Sigma _{t}^{2}\mathrm{d}t\right) \\
& \rightarrow \mathcal{MN}\left( 0,\frac{1}{3}\int_{0}^{1}\Gamma
_{t}^{4}\Sigma _{t}^{4}\mathrm{d}t\right) .
\end{split}%
\end{equation*}%
By Lemma~\ref{lem0}, it suffices to show 
\begin{equation*}
\sqrt{n}\sum_j\Gamma _{t_{j}}^{2}\Sigma _{t_{j}}^{2}\left(
\int_{t_{j}}^{t_{j+1}}n(W_{t}-W_{t_{j}})^{2}\mathrm{d}t-\frac{1}{2n}\right)
\rightarrow \mathcal{MN}\left( 0,\frac{1}{3}\int_{0}^{1}\Gamma
_{t}^{4}\Sigma _{t}^{4}\mathrm{d}t\right) .
\end{equation*}%
The left hand side is equal to 
\begin{equation*}
2n^{3/2}\sum_j\Gamma _{t_{j}}^{2}\Sigma
_{t_{j}}^{2}\int_{t_{j}}^{t_{j+1}}\int_{t_{j}}^{t}(W_{u}-W_{t_{j}})\mathrm{d}%
W_{u}\mathrm{d}t=:\sum_jA_{j}.
\end{equation*}%
Since 
\begin{equation*}
E_{t_{j}}[A_{j}]=E_{t_{j}}[A_{j}(W_{t_{j+1}}-W_{t_{j}})]=0
\end{equation*}%
and 
\begin{equation*}
E_{t_{j}}[A_{j}^{2}]=\frac{1}{3n}\Gamma _{t_{j}}^{4}\Sigma _{t_{j}}^{4},
\end{equation*}%
the result follows from Jacod's theorem of stable convergence. From this
computation, the uniform integrability is also clear.\newline

\noindent Step 2). Next, we show 
\begin{equation*}
n^{3/2}\sum_j\Gamma
_{t_{j}}^{2}\int_{t_{j}}^{t_{j+1}}(W_{t}-W_{t_{j}})^{2}(\Sigma
_{t}^{2}-\Sigma _{t_{j}}^{2})\mathrm{d}t
\end{equation*}%
is uniformly integrable and converges to $0$ in probability. The boundedness
in $L^2$ is not difficult to see, from which the uniform integrability
follows. Since 
\begin{equation*}
E\left[ \left\vert
\int_{t_{j}}^{t_{j+1}}(W_{t}-W_{t_{j}})^{2}\int_{t_{j}}^{t}\check{\Sigma}_{u}%
\mathrm{d}u\mathrm{d}t\right\vert \right] \leq n^{-3/2}\sqrt{ E\left[%
\int_{t_j}^{t_{j+1}}\left(\int_{t_j}^t \check{\Gamma}_u\mathrm{d}u\right)^2%
\mathrm{d}t\right] }
\end{equation*}%
by the Cauchy-Schwarz inequality, we have 
\begin{equation*}
E\left[\left| n^{3/2}\sum_j\Gamma
_{t_{j}}^{2}\int_{t_{j}}^{t_{j+1}}(W_{t}-W_{t_{j}})^{2}\int_{t_j}^t\check{%
\Sigma}_u\mathrm{d}u\mathrm{d}t\right|\right] \leq \frac{\|\check{\Sigma}\|_2%
}{\sqrt{n}} \to 0.
\end{equation*}
Therefore, it suffices to show 
\begin{equation*}
n^{3/2}\sum_j\Gamma
_{t_{j}}^{2}\int_{t_{j}}^{t_{j+1}}(W_{t}-W_{t_{j}})^{2}\int_{t_{j}}^{t}\hat{%
\Sigma}_{u}\mathrm{d}W_{u}\mathrm{d}t\rightarrow 0
\end{equation*}%
in probability. To see this, note that 
\begin{equation*}
\begin{split}
& E\left[ \left\vert n^{3/2}\sum_j\Gamma _{t_{j}}^{2}\hat{\Sigma}%
_{t_{j}}\int_{t_{j}}^{t_{j+1}}(W_{t}-W_{t_{j}})^{3}\mathrm{d}t\right\vert
^{2}\right] \\
& \leq Cn^{3}E\left[ \sum_j \Gamma _{t_{j}}^{4}\hat{\Sigma}_{t_{j}}^2
E_{t_{j}}\left[ \left( \int_{t_{j}}^{t_{j+1}}(W_{t}-W_{t_{j}})^{3}\mathrm{d}%
t\right) ^{2}\right] \right] \leq \frac{C}{n}
\end{split}%
\end{equation*}%
by  (H1) and (H3).
Therefore, with the aid of  (H3), it follows
from 
\begin{equation*}
\begin{split}
& E\left[ \left\vert \sum_j\Gamma
_{t_{j}}^{2}\int_{t_{j}}^{t_{j+1}}(W_{t}-W_{t_{j}})^{2}\int_{t_{j}}^{t}(\hat{
\Sigma}_{u}-\hat{\Sigma}_{t_{j}})\mathrm{d}W_{u}\mathrm{d}t\right\vert %
\right] \\
&= E\left[ \left| \int_0^T\Gamma_{[nt]/n}(W_t-W_{[nt]/n})^2\int_{[nt]/n}^t(%
\hat{\Sigma}_u-\hat{\Sigma}_{[nt]/n})\mathrm{d}W_u \mathrm{d}t \right| %
\right] \\
& \leq E\left[\int_0^T\Gamma_{[nt]/n}^4\mathrm{d}t\right]^{1/4} E\left[%
\int_0^T |W_t-W_{[nt]/n}|^8\mathrm{d}t\right]^{1/4} E\left[ \int_0^T
\int_{[nt]/n}^t |\hat{\Sigma}_u-\hat{\Sigma}_{[nt]/n}|^2\mathrm{d}u \right]%
^{1/2} \\
& \leq Cn^{-3/2}\sup_{0\leq |s-t|\leq 1/n}E[|\hat{\Sigma}_{s}-\hat{\Sigma}%
_{t}|^{2}]^{1/2} = o(n^{-3/2}).
\end{split}%
\end{equation*}

\noindent Step 3). Next, we look at 
\begin{equation*}
2n^{3/2}\sum_j\Gamma
_{t_{j}}\int_{t_{j}}^{t_{j+1}}(W_{t}-W_{t_{j}})\int_{t_{j}}^{t}(\Gamma
_{u}-\Gamma _{t_{j}})\mathrm{d}W_{u}\Sigma _{t}^{2}\mathrm{d}t.
\end{equation*}%
By a similar argument to the above, the problem reduces to showing that 
\begin{equation*}
2n^{3/2}\sum_j\Gamma _{t_{j}}\Sigma
_{t_{j}}^{2}\int_{t_{j}}^{t_{j+1}}(W_{t}-W_{t_{j}})\int_{t_{j}}^{t}%
\int_{t_{j}}^{u}\hat{\Gamma}_{s}\mathrm{d}W_{s}\mathrm{d}W_{u}\mathrm{d}t
\end{equation*}%
converges to $0$ in $L^{2}$, where $\hat{\Gamma}$ is the diffusion
coefficient of the continuous It\^{o} process $\Gamma $. Since 
\begin{equation*}
\begin{split}
& E_{t_{j}}\left[ \int_{t_{j}}^{t_{j+1}}(W_{t}-W_{t_{j}})\int_{t_{j}}^{t}%
\int_{t_{j}}^{u}\hat{\Gamma}_{s}\mathrm{d}W_{s}\mathrm{d}W_{u}\mathrm{d}t%
\right] \\
& =\int_{t_{j}}^{t_{j+1}}E_{t_{j}}\left[ \int_{t_{j}}^{t}\int_{t_{j}}^{u}%
\hat{\Gamma}_{s}\mathrm{d}W_{s}\mathrm{d}u\right] \mathrm{d}t=0,
\end{split}%
\end{equation*}%
we have 
\begin{equation*}
\begin{split}
& E\left[ \left\vert 2n^{3/2}\sum_j\Gamma _{t_{j}}\Sigma
_{t_{j}}^{2}\int_{t_{j}}^{t_{j+1}}(W_{t}-W_{t_{j}})\int_{t_{j}}^{t}%
\int_{t_{j}}^{u}\hat{\Gamma}_{s}\mathrm{d}W_{s}\mathrm{d}W_{u}\mathrm{d}%
t\right\vert ^{2}\right] \\
& \leq Cn^{3}\sum_jE\left[ \Gamma _{t_{j}}^2\Sigma_{t_{j}}^{4} \left\vert
\int_{t_{j}}^{t_{j+1}}(W_{t}-W_{t_{j}})\int_{t_{j}}^{t}\int_{t_{j}}^{u}\hat{%
\Gamma}_{s}\mathrm{d}W_{s}\mathrm{d}W_{u}\mathrm{d}t\right\vert ^{2}\right]
\\
& \leq Cn^{2}\sum_j\int_{t_{j}}^{t_{j+1}}E\left[ \Gamma
_{t_{j}}^2\Sigma_{t_{j}}^{4} (W_{t}-W_{t_{j}})^{2}\left(
\int_{t_{j}}^{t}\int_{t_{j}}^{u}\hat{\Gamma}_{s}\mathrm{d}W_{s}\mathrm{d}%
W_{u}\right) ^{2}\right] \mathrm{d}t \\
& \leq Cn^2 E\left[\int_0^T(\Gamma \Sigma^2)^8_{[nt]/n}\mathrm{d}t\right]%
^{1/4} E\left[ \int_0^T|W_t-W_{[nt]/n}|^8\mathrm{d}t \right]^{1/4} \\
& \hspace{2cm} \times E\left[ \int_0^T\left| \int_{[nt]/n}^t\int_{[nt]/t}^u%
\hat{\Gamma}_s\mathrm{d}W_s\mathrm{d}W_u \right|^4 \mathrm{d}t \right]^{1/2}
\\
& \leq Cn E\left[ \int_0^T\left| \int_{[nt]/n}^t \left|\int_{[nt]/t}^u\hat{%
\Gamma}_s\mathrm{d}W_s\right|^2\mathrm{d}u \right|^2 \mathrm{d}t \right]%
^{1/2} \\
& \leq Cn^{1/2} E\left[ \int_0^T \int_{[nt]/n}^t \left|\int_{[nt]/t}^u\hat{%
\Gamma}_s\mathrm{d}W_s\right|^4\mathrm{d}u \mathrm{d}t \right]^{1/2} \\
& \leq Cn^{1/2} E\left[ \int_0^T \int_{[nt]/n}^t \left|\int_{[nt]/t}^u\hat{%
\Gamma}_s^2\mathrm{d}s\right|^2\mathrm{d}u \mathrm{d}t \right]^{1/2} \leq 
\frac{C}{n}\|\hat{\Gamma}\|_4^{2} \to 0.
\end{split}%
\end{equation*}

\noindent Step 4). Next, we observe that 
\begin{equation*}
n^{3/2}\sum_j\int_{t_{j}}^{t_{j+1}}\left( \int_{t_{j}}^{t}(\Gamma
_{u}-\Gamma _{t_{j}})\mathrm{d}W_{u}\right) ^{2}\Sigma _{t}^{2}\mathrm{d}t
\end{equation*}%
is negligible. This simply follows from 
\begin{equation*}
E[|\Gamma _{u}-\Gamma _{t_{j}}|^{4}]\leq C|u-t_{j}|^2
\end{equation*}%
and so, we omit the detail. \newline

\noindent Step 5). It remains to show that the part involved with $\Xi$ ; 
\begin{equation*}
n^{3/2}\sum_j\int_{t_{j}}^{t_{j+1}}\left( \int_{t_j}^t\Xi_u\mathrm{d}%
u\right) ^{2}\Sigma _{t}^{2}\mathrm{d}t +
2n^{3/2}\sum_j\int_{t_{j}}^{t_{j+1}} \int_{t_j}^t\Xi_u\mathrm{d}u
\int_{t_j}^t\Gamma_u\mathrm{d}W_u\Sigma _{t}^{2}\mathrm{d}t
\end{equation*}
is negligible. The first term is easy to treat. For the second term, we
first observe that it can be approximated by 
\begin{equation*}
2n^{3/2}\sum_j\int_{t_{j}}^{t_{j+1}} \Xi_{t_j}(t-t_j)\Gamma_{t_j}
(W_t-W_{t_j})\Sigma_{t_j}^{2}\mathrm{d}t
\end{equation*}
by a similar argument as before. Then, using that 
\begin{equation*}
E_{t_j}\left[ \int_{t_{j}}^{t_{j+1}} \Xi_{t_j}(t-t_j)\Gamma_{t_j}
(W_t-W_{t_j})\Sigma_{t_j}^{2}\mathrm{d}t \right] = 0,
\end{equation*}
we have 
\begin{equation*}
\begin{split}
& 4n^{3} E\left[\left|\sum_j\int_{t_{j}}^{t_{j+1}}
\Xi_{t_j}(t-t_j)\Gamma_{t_j} (W_t-W_{t_j})\Sigma_{t_j}^{2}\mathrm{d}%
t\right|^2\right] \\
&=4n^{3} E\left[\sum_j \left|\int_{t_{j}}^{t_{j+1}}
\Xi_{t_j}(t-t_j)\Gamma_{t_j} (W_t-W_{t_j})\Sigma_{t_j}^{2}\mathrm{d}%
t\right|^2\right] \\
& \leq 4n^2 \sum_j E\left[\Xi_{t_j}^2\Gamma_{t_j}^2 \Sigma_{t_j}^{4}
\int_{t_{j}}^{t_{j+1}} (t-t_j)^2 (W_t-W_{t_j})^2\mathrm{d}t \right] \\
& = n^2\sum_j E[\Xi_{t_j}^2\Gamma_{t_j}^2 \Sigma_{t_j}^{4}](t_{j+1}-t_j)^4
\to 0
\end{split}%
\end{equation*}
since $t_j = t^n_j \wedge T$. \hfill////

\subsection{Limit results for $(V^n,Z^n)$}
\begin{lem}
\label{lemmaV}Consider $p>1$ and $n\geq 1$ and assume that (H1) holds. Then $V^n \in  \mathbb{L}_{-}^{1,p}$,\hspace{0.1cm} $D^{-}V^n\in\mathbb{L}^{1,p}$ and, for any  $\alpha>0$ and $\beta<1$ there exists a positive constant $C$ such that
\begin{equation}
\label{DnV}
 \|V^n\|_{1,p}<C n{^\alpha},
\end{equation}
\begin{equation}
\label{D-nV}
 \|D^{-}V^n\|_{1,p}<C n{^\alpha}
\end{equation}
and
\begin{equation}
\label{D-nVx}
 \|D_r V_t^n-D_r V_0^n\|_{p}<C n^{\alpha}(t-r)^{\beta}.
\end{equation}

\end{lem}
\textit{Proof:} We know that
$$V^n_t = n \int_t^T(X_s-X^n_s)^2\Sigma_s^2\mathrm{d}s.$$
Remark \ref{lem4p} gives us that $\Sigma^2 \in \mathbb{L}^{1,q}$ for every $q>1$. On the other hand, Proposition \ref{MalliavinderItoprocess}, Burkh\"older-Davis-Gundy inequality and H\"older's inequality give us that, for every $q>1$ and $\gamma>0$,  there exist a constant $C>0$ such that
\begin{equation}
\label{normofX}
\|X-X^n\|_{1,q}\leq C n^{\gamma -1/2}.
\end{equation}
Now (\ref{DnV}) results follows from a direct application of H\"older's inequaliy. A similar argument gives us (\ref{D-nV}). Finally, we can write
$$
D_rV_t^n-D_rV_0^n=-n\int_r^t D_r[(X_s-X_s^n)^2\Sigma_s^2]\mathrm{d}s,
$$
which, jointly with Burkh\"older-Davis-Gundy and H\"older inequalities, gives us (\ref{D-nVx}). Now the proof is complete.
\hfill////

\begin{lem}
\label{new} Assume that (H1) holds. Then, for every $n>1$, there exists an adapted and square integrable process $R^n$ satisfying that
\begin{equation}
\label{AB}
D^- V^n_s=R^n_s+n \int_s^T \left(\int_{\frac{[n\theta]}{n}}^\theta D_s(\Gamma_r^2 \Sigma_\theta^2)dr \right) \mathrm{d}\theta,
\end{equation}
 where,   for any $\delta<1$, there exist two constants $C$ and $p>1$ such that, for any $A\in \mathbb{D}^{2,p}$ and $s\in [0,T]$
\begin{equation}
\label{bond}
|E(AR_s^n)|\leq C n^{-\delta}\|A\|_{2,p}.
 \end{equation}
 
\end{lem}
\textit{Proof : } Note that
\begin{eqnarray*}
D^-V^n_s &=& 2n\int_s^T (X_\theta-X^n_\theta)[D_s(X_\theta-X^n_\theta)]\Sigma_\theta^2 \mathrm{d}\theta\\
&&+n \int_s^T (X_\theta-X^n_\theta)^2 D_s\Sigma_\theta^2 \mathrm{d}\theta .
\end{eqnarray*}
Now, as 
$$
Y^n_\theta:=X_\theta-X^n_\theta =\int_{\frac{[n\theta]}{n}}^\theta \Xi_\tau \mathrm{d}\tau +\int_{\frac{[n\theta]}{n}}^\theta  \Gamma_\tau \mathrm{d}W_\tau
$$
and 
$$
U^n_\theta:=D_s (X_\theta-X^n_\theta)=\bold{1}_{[\frac{[n\theta]}{n},\theta]} (s)\Gamma_s+\int_{\frac{[n\theta]}{n}}^\theta D_s \Xi_\tau \mathrm{d}\tau +\int_{\frac{[n\theta]}{n}}^\theta D_s \Gamma_\tau \mathrm{d}W_\tau,
$$
It\^{o}'s formula gives us that
\begin{eqnarray*}
&& D^-V^n_s= 2n\int_s^T Y_\theta^n U_\theta^n \Sigma_\theta^2 \mathrm{d}\theta+n\int_s^T (Y_\theta^n)^2 D_s\Sigma_\theta^2 \mathrm{d}\theta\\
&&=2n\int_s^T (X_\theta-X^n_\theta)  \bold{1}_{[\frac{[n\theta]}{n},\theta]} (s)\Gamma_s \Sigma_\theta^2 \mathrm{d}\theta\\
&&+2n \int_s^T\left(\int_{\frac{[n\theta]}{n}}^\theta \Xi_\tau \hat{U}_\tau^n \mathrm{d}\tau\right) \Sigma_\theta^2 \mathrm{d}\theta 
+2n \int_s^T \left(\int_{\frac{[n\theta]}{n}}^\theta \Gamma_\tau \hat{U}_\tau^n \mathrm{d}W_\tau\right) \Sigma_\theta^2 \mathrm{d}\theta\\
&&+2n\int_s^T \left(\int_{\frac{[n\theta]}{n}}^\theta D_s\Xi_\tau Y_\tau^n \mathrm{d}\tau\right) \Sigma_\theta^2 \mathrm{d}\theta+2n \int_s^T \left(\int_{\frac{[n\theta]}{n}}^\theta D_s\Gamma_\tau Y_\tau^n \mathrm{d}W_\tau\right) \Sigma_\theta^2 \mathrm{d}\theta\\
&&+2n \int_s^T\left(\int_{\frac{[n\theta]}{n}}^\theta \Xi_\tau Y_\tau^n \mathrm{d}\tau\right) D_s\Sigma_\theta^2 \mathrm{d}\theta +2n \int_s^T \left(\int_{\frac{[n\theta]}{n}}^\theta \Gamma_\tau Y_\tau^n \mathrm{d}W_\tau\right) D_s\Sigma_\theta^2 \mathrm{d}\theta\\
&&+2n \int_s^T \left(\int_{\frac{[n\theta]}{n}}^\theta \Gamma_\tau D_s\Gamma_\tau \mathrm{d}\tau\right) \Sigma_\theta^2 \mathrm{d}\theta+n \int_s^T \left(\int_{\frac{[n\theta]}{n}}^\theta \Gamma_\tau^2 \mathrm{d}\tau\right) D_s\Sigma_\theta^2\mathrm{d}\theta\\
&&=:R^n_s +2n \int_s^T \left(\int_{\frac{[n\theta]}{n}}^\theta \Gamma_\tau D_s\Gamma_\tau \mathrm{d}\tau\right) \Sigma_\theta^2 \mathrm{d}\theta+n \int_s^T \left(\int_{\frac{[n\theta]}{n}}^\theta \Gamma_\tau^2 \mathrm{d}\tau\right) D_s\Sigma_\theta^2\mathrm{d}\theta\\
&&=: R^n_s +n \int_s^T \left(\int_{\frac{[n\theta]}{n}}^\theta D_s(\Gamma_\tau^2\Sigma_\theta^2)\mathrm{d}\tau\right) \mathrm{d}\theta,
\end{eqnarray*}
where $\hat{U}^n_\tau=U^n_\tau- \bold{1}_{[\frac{[n\theta]}{n},\theta]} (s)\Gamma_s$.
Now, using again the duality relationship and H\"older's inequality  it is easy to see that
 $$ | E(A R_s^{n} )|\leq C n^{-1+\frac2q}\|A\|_{2,p}$$ 
for any $p>q>1$. Now,  taking $q>\frac{2}{1-\delta}$ the proof is complete. \hfill///
\begin{lem}
\label{dv} 
\begin{equation*}
\lim_{n\to \infty} n\int_{0}^{T}\int_{[nu]/n}^{u}E\left[
q_{1,1}(V_{0}^{n},0)(D^{-}V_{s}^{n}-D^{-}V_{s})\Theta _{u} \Gamma _{s}%
\right] \mathrm{d}s\mathrm{d}u 
  =0.
\end{equation*}

\end{lem}
\textit{Proof : } Lemma \ref{new} gives us that 
\begin{equation*}
\begin{split}
& n\int_{0}^{T}\int_{[nu]/n}^{u}E\left[
q_{1,1}(V_{0}^{n},0)(D^{-}V_{s}^{n}-D^{-}V_{s})\Theta _{u} \Gamma _{s}%
\right] \mathrm{d}s\mathrm{d}u \\
&= n\int_{0}^{T}\int_{[nu]/n}^{u}E\left[
q_{1,1}(V_{0}^{n},0)R^n_s\Theta _{u} \Gamma _{s}%
\right] \mathrm{d}s\mathrm{d}u\\
&+ n^2\int_{0}^{T}\int_{[nu]/n}^{u}E\left[
q_{1,1}(V_{0}^{n},0) \left[\left(   \int_s^T \int_{\frac{[n\theta]}{n}}^\theta D_s(\Gamma_\tau^2\Sigma_\theta^2)\mathrm{d}\tau \mathrm{d}\theta \right)-\frac{1}{n} D^-V_s\right]\Theta _{u} \Gamma _{s}%
\right] \mathrm{d}s\mathrm{d}u.
\end{split}
\end{equation*}
 (\ref{bond}) gives us that the first term in the right-hand side of the above equatlity tends to zero. On the other hand,  for any $s<\theta$, the process  $D_s(\Gamma^2\Sigma_\theta^2)$ is continuous in $L^p (\Omega)$. Then, a direct application of H\"{o}lder's inequality gives us that the second term tends to zero. Now the proof is complete.  \hfill////


\begin{lem}
\label{ddv} 
$$ n\int_{0}^{T}E[q_{1,1}(V_{0}^{n},0)(D^{-}V_{t}^{n}-D^{-}V_{t})
(X_{t}-X_{t}^{n})\Sigma _{t}]\mathrm{d}t\to \frac{1}{3}  E\left[q_{1,1}(V_{0},0)\int_0^T \Sigma_t^3\Gamma_t^3 \mathrm{d}t\right].
$$
\end{lem} 
\textit{Proof : }  We can make use of of the  computations in the proof of Lemma 5.5 to see that the leading terms should be
\begin{eqnarray*}
&&2n^2 \int_0^T \left(\int_{\frac{[ns]}{n}}^s E[q_{1,1}(V^n_0,0) (X_s-X^n_s)\Gamma_t \Sigma_s^2 (X_t-X^n_t)\Sigma_t]\mathrm{d}t\right)\mathrm{d}s\\
&&+n^2\int_0^T  E[q_{1,1}(V^n_0,0) \left(\int_t^T\int_{\frac{[ns]}{n}}^s D_t(\Gamma_\tau^2\Sigma_s^2)\mathrm{d}\tau\mathrm{d}s-\frac{1}{n} D^-V_t\right)  (X_t-X^n_t)\Sigma_t]\mathrm{d}t.
\end{eqnarray*}
By the duality relationship, the second term tends to zero. For the first one, we have that its limit is 
\begin{eqnarray*}
&&2 \lim_{n\to\infty} n^2 E \sum_{i=1}^n\int_{t_{i}}^{t_{i+1}} \int_{\frac{[ns]}{n}}^s E[q_{1,1}(V^n_0,0) (X_s-X^n_s)\Gamma_s \Sigma_s^2 (X_t-X^n_t)\Sigma_t]
\mathrm{d}t \mathrm{d}s \\
&& =2 \lim_{n\to\infty} n^2 E \sum_{i=1}^n\int_{t_{i}}^{t_{i+1}} \int_{\frac{[ns]}{n}}^s E[q_{1,1}(V^n_0,0) \Gamma_s \Sigma_s^2 \left(\int_{t_i}^t \Gamma_r^2 \mathrm{d}r\right)\Sigma_t]\mathrm{d}t \mathrm{d}s \\
&&=\frac{1}{3}\lim_{n\to\infty}  E\left[q_{1,1}(V_{0}^{n},0)\int_0^T \Sigma_t^3\Gamma_t^3 \mathrm{d}t\right],
\end{eqnarray*}
and this allows us to complete the proof. \hfill////

\vspace{0.3cm}
The proofs of  Lemmas \ref{lem2}  and \ref{lem3} will be based on the following technical result.
\begin{lem}
\label{elprocésZ} Suppose that (H1) holds. Consider a real function $f=q_{i,j}$, for some $i,j\geq 0$.
 Let $p > 1$. Then, for any $\alpha > 1/p$, there exists  $C > 0$  such that for any
 $A\in \mathbb{D}^{2,p}$,
$$
E\left[\left(f(V^n_t,Z^n_t)-f(V^n_0,0)\right)A\right]\leq C n^{\alpha-\frac{1}{2}}\|A\|_{2,p}.
$$
\end{lem}
\textit{Proof:}  Using that $f=q_{i,j}$ solves the heat equation, a direct application of the anticipating It\^{o}'s formula gives us that
\begin{eqnarray}
\label{MalliavinZ}
&&Af(V^n_t,Z^n_t)-Af(V^n_0,0)\nonumber \\
&&=\sqrt{n}\int_0^t A\partial_y f (V^n_s,Z^n_s)(X_s-X_s^n)\mathrm{d}Y_s\nonumber\\
&&+\sqrt{n}\int_0^t A\partial^2_{xy} f (V^n_s,Z^n_s)D^{-}V_s^n (X_s-X_s^n)\Sigma_s \mathrm{d}s\nonumber\\
&&+\sqrt{n}\int_0^t \partial_y f (V^n_s,Z^n_s)D_s A (X_s-X_s^n)\Sigma_s \mathrm{d}s.
\end{eqnarray}
Then, taking conditional expectations we get
\begin{eqnarray}
\label{MalliavinZexpectation}
&&E\left[\left(f(V^n_t,Z^n_t)-f(V^n_0,0)\right)A\right] \nonumber\\
&&=\sqrt{n}E\left[ \int_0^t A\partial_y f (V^n_s,Z^n_s)(X_s-X_s^n)\Theta_s\mathrm{d}s\right.\nonumber\\
&&+\int_0^t A\partial^2_{xy} f (V^n_s,Z^n_s)D^-V_s^n (X_s-X_s^n)\Sigma_s \mathrm{d}s+\left.\int_0^t \partial_y f (V^n_s,Z^n_s)D_s A (X_s-X_s^n)\Sigma_s \mathrm{d}s \right]\nonumber\\
&&=:\sqrt{n}\int_0^t  (X_s-X_s^n)J_s  \mathrm{d}s,
\end{eqnarray}
where
$$
J_s:=A\partial_y f (V^n_s,Z^n_s)\Theta_s+A\partial^2_{xy} f (V^n_s,Z^n_s)D^-V_s^n\Sigma_s+\partial_y f (V^n_s,Z^n_s)D_s A \Sigma_s.
$$
Now, the duality relationship between the Skorohod integral and the Malliavin derivative operator give us that
\begin{eqnarray}
&&E\left[\left(f(V^n_t,Z^n_t)-f(V^n_0,0)\right)A\right] \nonumber\\
&&=\sqrt{n}E\left[\int_0^t  \int_{\frac{[nt]}{n}}^t J_s \Xi_r \mathrm{d}s\mathrm{d}r+\int_0^t  \int_{\frac{[nt]}{n}}^t (D_r J_s) \Gamma_r \mathrm{d}s\mathrm{d}r\right]\nonumber\\
&&\leq n^{-\frac12+\alpha}E\left[\int_0^t  \int_{\frac{[nt]}{n}}^t |J_s \Xi_r|^{ \frac{1}{\alpha}} \mathrm{d}s\mathrm{d}r+\int_0^t  \int_{\frac{[nt]}{n}}^t |(D_r J_s) \Gamma_r |^{ \frac{1}{\alpha}}\mathrm{d}s\mathrm{d}r\right]^{\alpha}
\end{eqnarray}
for any $\alpha > 1/p$.
Then, the result follows as a direct consequence of (H1), Lemma \ref{lemmaV} and H\"{o}lder's inequality.
\hfill////

\section{The leading terms}
\label{leading}
Now we are in a position to prove the limit lemmas.
\begin{lem}
\label{lemT1} 
\begin{equation*}
\sqrt{n} \left\{ E[q(V^n_0,0)] - E[q(V_0,0)] \right\} \to 0
\end{equation*}
as $n\to \infty$.
\end{lem}
\textit{Proof : } By Taylor's formula, 
\begin{equation*}
\sqrt{n} \left\{ E[q(V^n_0,0)] - E[q(V_0,0)] \right\} = E[\sqrt{n}%
(V^n_0-V_0)\int_0^1q_{1,0}(V_0 + (V^n_0-V_0)s,0)\mathrm{d}s].
\end{equation*}
Since $q_{1,0}$ is bounded, the result follows from Proposition~\ref{stableconv}.
\hfill////

\vspace{0.2cm}

\begin{lem}
\label{lem2} 
\begin{equation*}
\begin{split}
& n \int_0^T E[q_{0,1}(V^n_t,Z^n_t)(X_t-X^n_t)\Theta_t]\mathrm{d}t \\
& \to \frac{1}{2} \int_0^TE[q_{0,1}(V_0,0)(\Xi_t\Theta_t +
D^+\Theta_t\Gamma_t)]\mathrm{d}t + \frac{1}{2}
\int_0^TE[q_{1,1}(V_0,0)D^{-}V_t\Theta_t\Gamma_t]\mathrm{d}t
\end{split}%
\end{equation*}
as $n\to \infty$.
\end{lem}
\textit{Proof:} We can write
\begin{equation*}
\begin{split}
& n \int_0^T E[q_{0,1}(V^n_t,Z^n_t)(X_t-X^n_t)\Theta_t]\mathrm{d}t \\
&=n \int_0^T E\left[q_{0,1}(V^n_t,Z^n_t)\left(\int_{\frac{[nt]}{n}}^t \Xi_r\mathrm{d}r\right)\Theta_t\right]\mathrm{d}t +n \int_0^T E\left[q_{0,1}(V^n_t,Z^n_t)\left(\int_{\frac{[nt]}{n}}^t \Gamma_r\mathrm{d}W_r\right)\Theta_t\right]\mathrm{d}t \\
&=n \int_0^T E\left[q_{0,1}(V^n_0,0)\left(\int_{\frac{[nt]}{n}}^t \Xi_r\mathrm{d}r\right)\Theta_t\right]\mathrm{d}t+n \int_0^T E\left[q_{0,1}(V^n_0,0)\left(\int_{\frac{[nt]}{n}}^t \Gamma_r\mathrm{d}W_r\right)\Theta_t\right]\mathrm{d}t\\
&+n \int_0^T E\left[(q_{0,1}(V^n_t,Z^n_t)-q_{0,1}(V^n_0,0))\left(\int_{\frac{[nt]}{n}}^t \Xi_r\mathrm{d}r\right)\Theta_t\right]\mathrm{d}t\\
&+n \int_0^T E\left[(q_{0,1}(V^n_t,Z^n_t)-q_{0,1}(V^n_0,0))\left(\int_{\frac{[nt]}{n}}^t \Gamma_r\mathrm{d}W_r\right)\Theta_t\right]\mathrm{d}t\\
&=:T_1+T_2+T_3+T_4.
\end{split}%
\end{equation*}
Now, 
\begin{equation*}
T_{1}\rightarrow \frac{1}{2}E\Big[q_{0,1}(V_{0},0)\int_{0}^{T}\Xi _{t}\Theta
_{t}\mathrm{d}t\Big]
\end{equation*}%
since 
\begin{equation*}
\begin{split}
& E\Big{|}\int_{0}^{T}n\int_{[nu]/n}^{u}[q_{0,1}(V_{0}^{n},0)\Xi _{s}\Theta
_{u}]\mathrm{d}s\mathrm{d}u-\frac{1}{2}\int_{0}^{T}[q_{0,1}(V_{0},0)\Xi
_{t}\Theta _{t}]\mathrm{d}t\Big{|} \\
& \leq E\Big{|}q_{0,1}(V_{0}^{n},0)\Big(\int_{0}^{T}n\int_{[nu]/n}^{u}\Theta
_{u}\Xi _{s}\mathrm{d}s\mathrm{d}u-\frac{1}{2}\int_{0}^{T}\Xi _{u}\Theta _{u}%
\mathrm{d}u\Big)\Big{|} \\
& +\frac{1}{2}E\Big{|}\Big(q_{0,1}(V_{0}^{n},0)-q_{0,1}(V_{0}^n,0)\Big)%
\int_{0}^{T}\Xi _{t}\Theta _{t}\mathrm{d}t\Big{|},
\end{split}%
\end{equation*}%
which tends to zero due to (H1), (H3), Lemma \ref{Hol} and Proposition 5.1.

For the second term, by the duality between the Malliavin derivative and the Skorohod integral,

\begin{equation*}
\begin{split}
&T_2=n\int_{0}^{T}\int_{[nu]/n}^{u}E[D_{s}[q_{0,1}(V_{0}^{n},0)\Theta
_{u}]\Gamma _{s}]\mathrm{d}s\mathrm{d}u \\
& =n\int_{0}^{T}\int_{[nu]/n}^{u}E\left[ \left(
q_{1,1}(V_{0}^{n},0)D^{-}V_{s}^{n}\Theta
_{u}+q_{0,1}(V_{0}^{n},0)D_{s}\Theta _{u}\right) \Gamma _{s}\right] \mathrm{d%
}s\mathrm{d}u \\
& =n\int_{0}^{T}\int_{[nu]/n}^{u}E\left[ \left(
q_{1,1}(V_{0}^{n},0)D^{-}V_{s}\Theta _{u}\right) \Gamma _{s}\right] ]\mathrm{%
d}s\mathrm{d}u \\
& +n\int_{0}^{T}\int_{[nu]/n}^{u}E\left[ \left(
q_{0,1}(V_{0}^{n},0)D_{s}^{+}\Theta _{u}\right) \Gamma _{s}\right] \mathrm{d}%
s\mathrm{d}u \\
& +n\int_{0}^{T}\int_{[nu]/n}^{u}E\left[ \left(
q_{1,1}(V_{0}^{n},0)(D^{-}V_{s}^{n}-D^{-}V_{s})\Theta _{u}\right) \Gamma _{s}%
\right] \mathrm{d}s\mathrm{d}u \\
& +n\int_{0}^{T}\int_{[nu]/n}^{u}E\left[ q_{0,1}(V_{0}^{n},0)\left(
D_{s}^{+}\Theta _{u}-D_{s}\Theta _{u}\right) \Gamma _{s}\right] \mathrm{d}s%
\mathrm{d}u
\end{split}%
\end{equation*}%
Notice that (H1), (H2) and   Lemma \ref{dv} imply that the last two terms  in the
above equality tend to zero. Then, similar arguments as for $T_{1}$ give us
that 
\begin{equation*}
T_{2}\rightarrow \frac{1}{2}\int_{0}^{T}E[q_{1,1}(V_{0},0)D^{-}V_{t}\Theta
_{t}\Gamma _{t}]\mathrm{d}t+\frac{1}{2}\int_{0}^{T}E[q_{0,1}(V_{0},0)D^{+}%
\Theta _{t}\Gamma _{t}]\mathrm{d}t.
\end{equation*}
Let us study $T_3$. We have that
\begin{eqnarray}
&&T_3= n \int_0^T E\left[(q_{0,1}(V^n_t,Z^n_t)-q_{0,1}(V^n_0,0))\left( \int_{\frac{[nt]}{n}}^t\Xi_r\mathrm{d}r\right)\Theta_t\right]\mathrm{d}t\nonumber\\
&&= n \int_0^T E\int_{\frac{[nt]}{n}}^t (q_{0,1}(V^n_t,Z^n_t)-q_{0,1}(V^n_0,0))\Xi_r\Theta_t\mathrm{d}r\mathrm{d}t,
\end{eqnarray}
which tends to zero due to  Lemma \ref{elprocésZ}.
For the last term,
\begin{eqnarray}
&&T_4=n \int_0^T E\left[(q_{0,1}(V^n_t,Z^n_t)-q_{0,1}(V^n_0,0))\left(\int_{\frac{[nt]}{n}}^t \Gamma_r\mathrm{d}W_r\right)\Theta_t\right]\mathrm{d}t \nonumber\\
&&=n \int_0^T E\left[\int_{\frac{[nt]}{n}}^t D_r\left[\Theta_t(q_{0,1}(V^n_t,Z^n_t)-q_{0,1}(V^n_0,0))\right] \Gamma_r\mathrm{d}r\right]\mathrm{d}t\nonumber\\
&&=n \int_0^T E\left[\int_{\frac{[nt]}{n}}^t \left(D_r\Theta_t\right) (q_{0,1}(V^n_t,Z^n_t)-q_{0,1}(V^n_0,0)) \Gamma_r\mathrm{d}r\right]\mathrm{d}t\nonumber\\
&&+n \int_0^T E\left[\int_{\frac{[nt]}{n}}^t [\Theta_t D_r\left(q_{0,1}(V^n_t,Z^n_t)-q_{0,1}(V^n_0,0))\right] \Gamma_r\mathrm{d}r\right]\mathrm{d}t\nonumber\\
&&=:T_4^1+T_4^2.
\end{eqnarray}
Using again (H1) and  Lemma \ref{elprocésZ}, we can easily check that $T_4^1\rightarrow 0$.  Now, 
\begin{eqnarray}
&&T_4^2=n \int_0^T E\left[\int_{\frac{[nt]}{n}}^t \Theta_t  \left(q_{1,1}(V^n_t,Z^n_t) D_rV_t^n-q_{1,1}(V^n_0,0) D_rV_0^n\right)\Gamma_r\mathrm{d}r\right]\mathrm{d}t\nonumber\\
&&+n \int_0^T E\left[\int_{\frac{[nt]}{n}}^t \Theta_t  \left(q_{0,2}(V^n_t,Z^n_t)D_rZ^n_t\right) \Gamma_r\mathrm{d}r\right]\mathrm{d}t\nonumber\\
&&=n \int_0^T E\left[\int_{\frac{[nt]}{n}}^t \Theta_t  \left(q_{1,1}(V^n_t,Z^n_t) -q_{1,1}(V^n_0,0) \right)D_rV_t^n\Gamma_r\mathrm{d}r\right]\mathrm{d}t\nonumber\\
&& +n \int_0^T E\left[\int_{\frac{[nt]}{n}}^t \Theta_t  q_{1,1}(V^n_0,0)\left(D_rV_t^n-D_rV_0^n\right)\Gamma_r\mathrm{d}r\right]\mathrm{d}t\nonumber \\
&&+n \int_0^T E\left[\int_{\frac{[nt]}{n}}^t \Theta_t  \left(q_{0,2}(V^n_t,Z^n_t)D_rZ^n_t\right) \Gamma_r\mathrm{d}r\right]\mathrm{d}t\nonumber\\
&&=:T_4^{2,1}+T_4^{2,2}+T_4^{2,3}.
\end{eqnarray}
The same arguments as before, together with Lemma \ref{lemmaV} gives us that $T_4^{2,1}+T_4^{2,2}\rightarrow 0$. On the other hand,

\begin{eqnarray}
&& D_rZ_t^n \nonumber \\
&=&\sqrt{n} D_r \left[\int_0^t (X_s-X_s^n)\Theta_s \mathrm{d}s+\int_0^t (X_s-X_s^n)\Sigma_s \mathrm{d}W_s\right]\nonumber \\
&=&\sqrt{n}\left[\int_r^t D_r(X_s-X_s^n)\Theta_s \mathrm{d}s+\int_r^t (X_s-X_s^n)D_r\Theta_s \mathrm{d}s\right.\nonumber \\
&+&\left. (X_r-X_r^n)\Sigma_r +\int_r^t D_r(X_s-X_s^n)\Sigma_s \mathrm{d}W_s+\int_r^t (X_s-X_s^n)D_r\Sigma_s \mathrm{d}W_s\right]\nonumber \\
&=:& (X_r-X_r^n)\Sigma_r +H_{r,t}.
\end{eqnarray}
Then
\begin{eqnarray}
\label{laultima}
&&T_4^{2,3}=n \int_0^T E\left[\int_{\frac{[nt]}{n}}^t \Theta_t  \left(q_{0,2}(V^n_t,Z^n_t)(X_r-X_r^n)\Sigma_r \right) \Gamma_r\mathrm{d}r\right]\mathrm{d}t\nonumber\\
&&+n \int_0^T E\left[\int_{\frac{[nt]}{n}}^t \Theta_t  \left(q_{0,2}(V^n_t,Z^n_t) H_{r,t}\Sigma_r \right) \Gamma_r\mathrm{d}r\right]\mathrm{d}t.
\end{eqnarray}
Using again the duality relationship we deduce that the first term in the right-hand side of (\ref{laultima}) tends to zero. Moreover, Bulkh\"{o}lderDavis-Gundy and H\"{o}lder inequalities give us that for any $p>1, \gamma<1/2$,  $\| H_{r,t}\|_{p}\leq L_{r} (t-r)^{\gamma}$, for some process $L\in L^p$. This allows us to complete the proof.
\hfill ////

\begin{lem}
\label{lem3} 
\begin{equation}
\label{limlem3}
\begin{split}
& n \int_0^T E[q_{1,1}(V^n_t,Z^n_t)D^{-}V^n_t(X_t-X^n_t)\Sigma_t]\mathrm{d}t \\
& \to \frac{1}{2}\int_0^TE[q_{1,1}(V_0,0) (D^{-}V_t\Xi_t \Sigma_t +
(D^{-})^2V_t\Sigma_t \Gamma_t + D^{-}_tV_tD^{+}\Sigma_t\Gamma_t) ]\mathrm{d}t \\
& \hspace*{2cm} + \frac{1}{2} \int_0^TE[q_{2,1}(V_0,0)|D^{-}V_t|^2\Gamma_t%
\Sigma_t]\mathrm{d}t + \frac{1}{3} \int_0^T E[q_{1,1} (V_0,0) \Sigma_t^3\Gamma_t^3 ]\mathrm{d}t
\end{split}%
\end{equation}
as $n\to\infty$.
\end{lem}
\textit{Proof : }
We can write
\begin{eqnarray}
&& n \int_0^T E[q_{1,1}(V^n_t,Z^n_t)D^{-}V^n_t(X_t-X^n_t)\Sigma_t]\mathrm{d}t \nonumber \\
&&=n \int_0^T E[q_{1,1}(V^n_t,Z^n_t)D^{-}V_t(X_t-X^n_t)\Sigma_t]\mathrm{d}t \nonumber \\
&&+n \int_0^T E[q_{1,1}(V^n_t,Z^n_t)(D^{-}V^n_t-D^{-}V_t)(X_t-X^n_t)\Sigma_t]\mathrm{d}t \nonumber \\
&&=:T_1+T_2.
\end{eqnarray}
Notice that $T_1$ is similar to the term studied in Lemma \ref{lem2}, replacing $q_{0,1}$ by $q_{1,1}$ and $\Theta_t$ by 
$D^{-}V_t\Sigma_t$. Then, the same arguments as in the proof of the previous result give us that 
\begin{eqnarray}
T_1 &\to & \frac{1}{2} \int_0^TE[q_{1,1}(V_0,0)(\Xi_t(D^{-}V_t\Sigma_t) +
D^{-}(D^{-}V_t\Sigma_t)\Gamma_t)]\mathrm{d}t \nonumber\\
&&+ \frac{1}{2}
\int_0^TE[q_{2,1}(V_0,0)D^{-}V_t(D^{-}V_t\Sigma_t)\Gamma_t]\mathrm{d}t.
\end{eqnarray}
This coincides with  the first three terms in (\ref{limlem3}). For the second term
\begin{eqnarray}
T_2 &&=:n \int_0^T E[q_{1,1}(V^n_0,0)(D^{-}V^n_t-D^{-}V_t)(X_t-X^n_t)\Sigma_t]\mathrm{d}t \nonumber \\
&&+n \int_0^T E[(q_{1,1}(V^n_t,Z_t^n))-q_{1,1}(V^n_0,0)(D^{-}V^n_t-D^{-}V_t)(X_t-X^n_t)\Sigma_t]\mathrm{d}t \nonumber \\
&&=:T_2^1+T_2^2.
\end{eqnarray}
Lemma \ref{ddv} gives us that $T_2^1\to \frac{1}{3} E\left[ q_{1,1} (V_0,0) \int_0^T\Sigma_t^3\Gamma_t^3 dt\right]$. The same arguments as in the proof of this Lemma  allow us to prove that $T_2^2$ tends to zero. Now the proof is complete.\hfill ////

\end{document}